\documentclass[12pt,a4paper, twoside, reqno]{amsart} 

\usepackage[margin=3cm]{geometry}
\usepackage{ifpdf}
\usepackage{graphicx}
\usepackage[T1]{fontenc} 
\usepackage[latin1]{inputenc} 
\usepackage{amsmath} 
\usepackage{amsthm} 
\usepackage{amsfonts} 
\usepackage{amssymb}
\usepackage{amsbsy}
\usepackage{bbm} 
\usepackage[english]{babel} 
\usepackage{varioref} 
\usepackage{enumerate} 

\theoremstyle{definition}

\newtheorem*{defn*}{Definition}

\theoremstyle{remark}

\theoremstyle{plain}
\newtheorem{thm}{Theorem}
\newtheorem*{thm*}{Theorem}

\newtheorem*{thmktv}{Theorem KTV}

\newcommand{\Bad}{\text{\textup{\bf Bad}}}

\newcommand{\norm}[1]{\ensuremath{\left\Vert #1 \right\Vert}}

\newcommand{\infabs}[1]{\ensuremath{\left\vert #1 \right\vert}}

\newcommand{\hdim}{\text{\textup{dim}}_{\text{\textup{H}}}}

 \newcommand{\R}{\mathbb{R}}
  
   \newcommand{\RNI}{\mathbb{R}^{N+1}}
 \newcommand{\N}{\mathbb{N}}
 \newcommand{\Q}{\mathbb{Q}}
 
  \newcommand{\x}{\mathbf{x}}

            \newcommand{\qie}{\mathbf{q}_i}
            \newcommand{\qIe}{\mathbf{q}_1}

            \newcommand{\qNIe}{\mathbf{q}_{N+1}}
   \newcommand{\xp}{\mathbf{\hat{x}}}
    \newcommand{\qp}{\mathbf{\hat{q}}}
   \newcommand{\yp}{\mathbf{\hat{y}}}

          \newcommand{\qip}{\mathbf{\hat{q}}_i}
          \newcommand{\qIp}{\mathbf{\hat{q}}_1}

          \newcommand{\qNIp}{\mathbf{\hat{q}}_{N+1}}
   \newcommand{\PNR}{\mathbb{P}^N(\R)}
    \newcommand{\PNQ}{\mathbb{P}^N(\Q)}
   
   \newcommand{\zero}{\mathbf{0}}
     \newcommand{\Proj}{\mathbb{P}}

\title{A note on badly approximabe sets in projective space}

\author{Stephen Harrap$^\dag$}
\thanks{$^\dag$ The research is supported by EPSRC grant number $EP/L005204/1$ and the Danish research council.}

\address{S. Harrap,Durham University, Department of Mathematical Sciences, Science Laboratories, South Rd, Durham,
DH1 3LE, United Kingdom}
\email{s.g.harrap@durham.ac.uk}
\author{Mumtaz Hussain$^*$ }
\thanks{ $^*$ The research is supported by  the Australian research council}

\address{M. Hussain,School of Mathematical and Physical Sciences, The University of Newcastle,
Callaghan, NSW 2308, Australia.}
\email{mumtaz.hussain@newcastle.edu.au}

\subjclass[2000]{37A45, 11K60, 11J83, 11J86}

\begin{document}

\begin{abstract} 
	Recently, Ghosh \& Haynes \cite{HG} proved a Khintchine-type result for the problem of Diophantine approximation in certain projective spaces. In this note we complement their result by observing that a Jarn\'{\i}k-type result also holds for `badly approximable' points in real projective space. In particular, we prove that the natural analogue in projective space of the classical set of badly approximable numbers has full Hausdorff dimension when intersected with certain compact subsets of real projective space. Furthermore, we also establish an analogue of Khintchine's theorem for convergence relating to `intrinsic' approximation of points in these compact sets.

\end{abstract}

\maketitle

\section{Introduction}

It is an elementary fact that the rational numbers are dense as a subset of the real numbers. The fundamental aim of Diophantine approximation is to quantify this density by assigning an `arithmetical complexity' to each rational number; in the classical setup this complexity is simply given by the absolute value of denominator of the rational. The basic problem then amounts to choosing an irrational number~$x$ and asking in general at what `rate' we can find close approximations to $x$ by rationals of increasing complexity. Typically this `rate' is realized by some decreasing arithmetic function, as exemplified by the first classical result in the field.
\begin{thm}[Dirichlet, $1842$]
	\label{cor:dircor}
	For any real number $x$ there exist infinitely many rationals $p/q$ such that
	\begin{equation*}
	\label{eqn:dirichlet}
	\left|x-\frac{p}{q}\right|\, < \,\frac{c}{\infabs{q}^2}
	\end{equation*}
	when $c=1$.
\end{thm}
This famous result tells us that every for every irrational can be approximated by rationals at a rate defined by the reciprocal of the square of the denominator of each rational. Dirichlet's theorem was famously improved by Hurwitz in $1891$, who found an optimum choice for the constant $c$.
\begin{thm}[Hurwitz, $1891$]
	We may take $c=1/\sqrt{5}$ in the statement of Dirichlet's Theorem, and this choice of constant is best possible.
\end{thm}
The choice $1/\sqrt{5}$ is best possible in the following sense. For $c>0$ let $\Bad(c)$ denote the set of real numbers such that for all rationals $p/q$ we have
\begin{equation*}
\infabs{x- \frac{p}{q}} \, \geq \, \frac{c}{\infabs{q}^2}.
\end{equation*}
Then, for every $c$ satisfying $1/\sqrt{5}\geq c>0$ the set $\Bad(c)$ is nonempty. However, numbers of this type are rare in the sense that the set $\Bad:=\bigcup_{c>0}\Bad(c)$ of \textit{badly approximable numbers} has zero Lebesgue measure $\lambda$. One can show this directly using the famous Borel-Cantelli Lemma from probability theory, but it also follows from the seminal work of Khintchine, who gave a comprehensive description of the `typical' rates at which real numbers can be approximated by rationals.
\begin{thm}[Khintchine, 1926]
	Let $\psi: \N \rightarrow \R_{>0}$ be a non-increasing arithmetic function. Then, for $\lambda$-almost every real number $x$ the inequality
	$$
	\left|x-\frac{p}{q}\right|\, < \,\psi(\infabs{q})
	$$
	has infinitely many solutions $\frac{p}{q}\in \Q$ if and only if the sum $\sum_{m\in \N} \, m\, \psi(m)$ diverges.
\end{thm}
To complete out recollection of the classical results in metric Diophantine approximation, it was subsequently shown by Jarn\'{\i}k in $1928$ that the set $\Bad$ is in some sense as large as it can be; that is, in terms of its Hausdorff dimension, $\dim_H$.
\begin{thm}[Jarn\'{\i}k, 1928]
	The set $\Bad$ has full Hausdorff dimension; i.e., we have $\dim_H(\Bad)=1$.
\end{thm}
\noindent
Proofs of all these results, and various higher dimensional analogues, can be found in most classical texts on Diophantine approximation, for example \cite{Harm}.


 The concepts described above have an attractive rendition in the context of projective geometry. For example, consider the one-dimensional space of all straight lines in $\R^2$ passing through the origin. One may ask at what `rate' lines with rational slopes can approximate those with an irrational slopes. As before, to get `close' approximations we will require the lines with rational slope to be in some sense `complex'.  As it turns out, a suitable method of defining `closeness' here is in terms of the \emph{sine} of the (acute positive) angle between the lines in question. To describe how the classical results translate to this setting we first set the problem up precisely in the more general context of projective spaces over number fields.

 Let $k$ be a number field and let $k_v$ be its completion at the place $v$. Furthermore, let
 $\norm{ \cdot }$ denote an absolute value from the place $v$ that extends the Euclidean absolute value on $k_v$ if $v | \infty$ and the $p$-adic absolute value on $k_v$ if $v | p$. For notational convenience we work with  the normalized absolute value $\infabs{ \cdot }_v$ defined by
 $$
 \infabs{ \cdot }_v: \, = \, \norm{ \cdot }^{d_v/d},
 $$
 where $d=\left[ k \, : \, \Q\right]$ and $d_v=\left[ k_v \, : \, \Q_v\right].
 $
 This absolute value extends to a norm on $(N+1)$-dimensional vector spaces over $k_v$ in the following way. For any $\x=(x_0, \ldots, x_N)$ in $k_v^{N+1}$ let
 $$
 \infabs{\x}_v \, = \, \left( \sum_{j=0}^{N} \infabs{\x}_v^{2d/d_v} \right)^{\frac{d_v}{2d}}
 $$
 if $v | \infty$, and if $v | p$ let
  $$
  \infabs{\x}_v \, = \, \left( \max_{0 \leq j \leq N} \infabs{\x}_v^{d/d_v} \right)^{\frac{d_v}{d}}.
  $$
    The standard $N$-dimensional projective space over $k_v^{N+1}$, denoted by $\Proj^N(k_v)$, then comes equipped with the following bounded metric, as described in \cite{Rum}. Let $\phi: k_v^{N+1}\setminus \left\{\mathbf{0} \right\} \rightarrow \Proj^N(k_v)$ denote the standard quotient map; that is, $\phi$ assigns to each non-zero vector  $\x=(x_0, \ldots, x_N)$ the unique point $\xp \in \Proj^N(k_v)$ with homogeneous coordinates $\left[x_0, \ldots, x_N\right]$.  Then, for any given $\xp, \yp \in \Proj^N(\R)$ with $\xp= \phi(\mathbf{x})$ and $\yp= \phi(\mathbf{y})$ we define
$$
\delta_v(\xp, \yp):= \frac{\infabs{\mathbf{x} \wedge \mathbf{y}}_v}{\infabs{\mathbf{x}}_v \infabs{\mathbf{y}}_v},
$$
where $\wedge$ denotes the exterior (or wedge) product. One can check that this metric induces the usual quotient
topology on $\Proj^N(k_v)$.

As previously mentioned, in the case when $k=\Q$ and $k_v=\R$ the metric $\delta_\infty$ has a nice geometric interpretation in terms of the smallest (plane) angle between the lines spanned by the points $\mathbf{x}$ and $\mathbf{y}$ respectively -- it is given by the absolute value of the \emph{sine} of this angle.  Furthermore, viewing $\Proj^N(\R)$ as the set of all one-dimensional subspaces of $\R^{N+1}$, the map $\phi$ assigns to each point $\x \in \R^{N+1}\setminus \left\{\mathbf{0} \right\}$ the unique line in $ \R^{N+1}$ containing $\x$ and passing through the origin. It is easy to verify that a point $\xp \in \Proj^1(\R)$ is contained in the subspace $\Proj^1(\Q)$ if and only if $\delta_\infty(\xp, \phi((1,0)))$ takes a value in $\Q \cap [0,1]$; that is, iff the sine of the angle between the line spanning $\x$ and the line $y=0$ is rational. Similar characterisations hold when $N > 1$.

To discuss the process of approximation of points in $\Proj^N(k_v)$ by points in $\Proj^N(k)$ we will assign a notion of arithmetical complexity to each member of $\Proj^N(k)$. The following standard height function provides such a notion. For $\xp$ in $\Proj^N(k)$ let
$$
H(\xp): \, = \, \prod_{v} \, \infabs{\x}_v,
$$
where the product is taken over all places $v$. For $\qp  \in \Proj^N(\Q)$ with homogeneous coordinates $[q_1, \cdots, q_{N+1}]$ this reduces to
$$
H(\qp) = \max(\infabs{q'_1}, \ldots, \infabs{q'_{N+1}}),
$$
where $q'_1, \ldots, q'_{N+1}$ are integers (not all zero) satisfying $\gcd(q'_1, \ldots, q'_{N+1})=1$ and that for some non-zero real number $\alpha$ we have $q'_i= \alpha q_i$ for every $i=1 \ldots, N+1$.

The first results concerning Diophantine approximation in projective space were proven by Schmidt in the $1960$s, who established analogues of Dirichlet's Theorem and Hurwitz' Theorem in $\Proj^N(\R)$ (see \cite{Schmidt} Theorem $15$ Corollary and Theorem $18$ respectively). In particular, he proved the following statements describing the rate at which points in $\Proj^N(\R)$ can be approximated by `rational' points in $\Proj^N(\R)$; that is, points in the subspace $\Proj^N(\Q)$.
\begin{thm}[Schmidt \cite{Schmidt}, $1967$]
Fix $N\geq 1$. Then, there is a strictly positive absolute constant $c_N$ such that for every  $\xp \in  \Proj^N(\R) \setminus \Proj^N(\Q)$ there exist infinitely many rational points $\qp \in \Proj^N(\Q)$ satisfying
$$
\delta_\infty(\xp, \qp) \: < \: \frac{c_N}{H(\qp)^{(N+1)/N}}.
$$
\end{thm}
\begin{thm}[Schmidt \cite{Schmidt}, $1967$]
	Let $\alpha = (\sqrt{5}-1)/2$, and let $\xp_1,\dots, \xp_4 \in  \Proj^1(\R)$ be the points associated with the lines $x_1-\alpha x_2 = 0$, $\alpha x_1-x_2 = 0$, $x_1 + \alpha x_2 = 0$, $\alpha x_1 + x_2 = 0$ in $\R^2$, respectively. Given $\xp \in  \Proj^1(\R) \setminus \Proj^1(\Q)$ not among $\xp_1,\dots, \xp_4$, then there are infinitely many rational $\qp \in \Proj^1(\Q)$ with
	$$
	\delta_\infty(\xp, \qp) \: < \: \frac{1}{\sqrt{5} \, H(\qp)^2}.
	$$
	 Here the constant $1/\sqrt{5}$ is best possible. If $\xp$ is among $\xp_1,\dots, \xp_4$ then for any $\epsilon > 0$ there are infinitely many rational $\qp \in \Proj^1(\Q)$  satisfying
	 $$
	  \delta_\infty(\xp, \qp) \: < \: \frac{1}{\sqrt{5} \, H(\qp)^2} +  \frac{1/\sqrt{500}+\epsilon}{H(\qp)^6}.
	  $$
	  Here, $1/\sqrt{5}$ and $1/\sqrt{500}$ are best possible.
\end{thm}

In $1999$, Choi \& Vaaler extended the first of Schmidt's above results to the general number field setting we described earlier.
\begin{thm}[Choi \& Vaaler \cite{CV}, 1999]
	\label{thm:CV}
	Given a number field $k$, fix $N\geq 1$ and a (finite or infinite) place $v$ of $k$. Then, there is a strictly positive absolute constant $c_{N,v}$ such that for every  $\xp \in  \Proj^N(k_v) \setminus \Proj^N(k)$ there exist infinitely many rational points $\qp \in \Proj^N(k)$ satisfying
	$$
	\delta_v(\xp, \qp) \: < \: \frac{c_{N,v}}{H(\qp)^{(N+1)/N}}.
	$$
\end{thm}

Inspired by this result, in a recent paper \cite{HG} Haynes \& Ghosh established an analogue of Khintchine's theorem for projective spaces in the case $k=\Q$. Before stating their result, we must define a suitable probability measure $\hat\mu_v$ on $\Proj^N(k_v)$. Following the construction defined by Choi \cite{Choi}, as used in \cite{HG}, we first specify a measure $\mu_v \,(:=\mu_v^{N+1})$ on $k_v^{N+1}$. 

If $v$ is an infinite place then $\mu_v$ is taken to be the usual $(N+1)$-fold Lebesgue measure, and if $v$ is a finite place then we take $\mu_v$ to be the $(N+1)$-fold Haar measure, normalised so that
$$\mu_v^1(\mathcal{O}_v) \, = \, \infabs{\mathcal{D}_v}_v^{d/2}.
$$
Here, $\mathcal{O}_v$ is the ring of integers of $k_v$ and $\mathcal{D}_v$ is the local different of $k$ at $v$. We then simply define the $\sigma$-algebra $\mathcal{M}$ of sets in $\Proj^N(k_v)$ to be the collection of sets $M \subseteq \Proj^N(k_v)$ such that $\phi^{-1}(M)$ lies in the Borel $\sigma$-algebra of $k_v^{N+1}$. The measure $\hat\mu_v$ on $\Proj^N(k_v)$ is then defined for any $M \in \mathcal{M}$ by
\begin{equation}\label{eqn:projmeasure}
\hat\mu_v(M): \, = \, \frac{\mu_v( \phi^{-1}(M) \cap  B(\mathbf{0}, 1) ) }{\mu_v(B(\mathbf{0}, 1) )},
\end{equation}
where $B(\mathbf{0}, 1)$ is the standard closed unit ball in $k_v^{N+1}$, centred at the origin.
\begin{thm}[Haynes \& Ghosh \cite{HG}, 2013]
	Let $\psi: \N \rightarrow \R_{>0}$ be a non-increasing arithmetic function and let $v$ be a (finite or infinite) place of $\Q$. Then, for $\hat\mu_v$-almost every $\xp \in \Proj^N(\Q_v)$ the inequality
	$$
	\delta_v(\xp, \qp )\, < \,\psi(H(\qp))
	$$
	has infinitely many solutions $\qp\in \PNQ$ if and only if the sum $\sum_{m\in \N} \, m^N\, \psi(m)^N$ diverges.
\end{thm}



An important corollary of this result is that the exponent $-(N+1)/N$ in Theorem~\ref{thm:CV} is generically optimal when $k = \Q$. Furthermore, it implies that the natural analogue of the set of badly approximable numbers within $\Proj^N(\Q_v)$ has zero $\hat\mu_v$-measure; that is, we have
$\hat\mu_v(\Bad_\Proj^v)=0$, where
\begin{equation*}
\Bad_\Proj^v: \, = \Bad_\Proj^v(N): \, = \, \left\{\xp\in\Proj^N(\Q_v): \: \inf_{\qp \, \in \,  \Proj^N(\Q)} H(\qp)^{\frac{N+1}{N}} \, \delta_v(\xp, \qp)\, >\,  0 \right\}.
  \end{equation*}

  A special case of the main result of this note implies that (when $v=\infty$ and $\Q_v=\R$) the set $\Bad_\Proj^\infty$ has full Hausdorff dimension, and so the theory of Diophantine approximation in projective space in this setting fits in line with the classical results.
\begin{thm} For every $N\in\N$ the set $\Bad_\Proj^\infty(N)$ has full Hausdorff dimension $N$.
\end{thm}

In fact, our main result is deeper and extends to more general measures than those described above. Let $\Omega$ be a compact subset of $\Proj^N(\R)$ which supports a `well-behaved' measure (see later section and the theorem below for a precise conditions the measure must satisfy). Then we are able to show that the intersection of $\Bad_\Proj^\infty$ with $\Omega$ has full Hausdorff dimension inside $\Omega$.

\begin{thm} \label{thm:omegadim}
Let $\Omega$ be a compact subset of $\Proj^N(\R)$ which supports a $\delta$-Ahlfors
regular and projectively absolutely $\eta$-decaying measure $m$, for some $\delta, \eta>0$. Then,
$$
\hdim(\Bad_\Proj^\infty \cap \Omega) \: = \: \hdim(\Omega) \: = \: \delta.
$$
\end{thm}

To complement this result we also exhibit an analogue of the `convergence' part of Khintchine's theorem for approximation inside $\Omega$. In particular we demonstrate that the following projective analogue of a theorem of Weiss \cite{Barak} holds.
\begin{thm}\label{thm:Barak}
	Let $\Omega$ be a compact subset of $\Proj^N(\R)$ which supports a $\delta$-Ahlfors
	regular and projectively absolutely $\eta$-decaying measure $\hat{\mu}$, for some $\delta, \eta>0$, and let $\psi: \N \rightarrow \R_{>0}$ be a non-increasing arithmetic function. Then for $\hat{\mu}$-almost every  $\xp$ in $\Proj^N(\R) \cap \Omega$ the inequality
	$$
	\delta_\infty(\xp, \qp )\, < \,\psi(H(\qp))
	$$
	has only finitely many solutions if the sum
		$$
		\sum_{m=1}^\infty \, m^{\frac{(N+1)}{N}\eta-1}	\psi(m)^\eta$$
		converges.
\end{thm}
One consequence of this statement is that the set $\Bad_\Proj^\infty$ has zero $\hat{\mu}$-measure inside $\Omega$, and so Theorem \ref{thm:omegadim} is indeed non-trivial. The proof closely follows the method outlined in \cite{Barak} (see also \cite{PVafm}) and for the sake of completion we choose to include all the details. Unfortunately, a proof of the `divergent' counterpart to this result remains out of reach, a situation which is in total analogue with the classical setting.


\section{Measure Conditions}
\subsection{Definitions}\label{sec:measuredefs}

In their 2004 paper \cite{MR2134453}, Kleinbock, Lindenstrauss \& Weiss introduced the notion of `friendly' measures on compact subsets of $\R^{N+1}$, a class of measures adhering to certain rigid geometric conditions. A more restrictive subclass was
investigated by Pollington \& Velani in~\cite{PVafm} and it is with the equivalent ideas in projective space that we are concerned.

The following properties will be enforced on
any locally finite Borel measure $m$ supported on a compact subset $\Omega$ of $\Proj^N(\R)$.
Throughout, $\hat{B}(\xp, r)$ will
denote the closed ball in $\Proj^N(\R)$ (with respect to the metric $\delta_\infty$) with centre $\xp \in \Omega$ and radius $r>0$. For technical reasons we will only work with balls of radius $r \leq \sin1$.

Firstly, a measure $m$ is said to be \emph{$\delta$-Ahlfors
	regular} if there exist strictly positive constants $\delta$ and
$r_0$ such that for $\xp \in \Omega$ and $r< r_0$ we have
\begin{equation*}
	a r^{\delta} \leq m(\hat{B}(\xp, r)) \leq b r^{\delta},
\end{equation*}
where $0 < a \leq 1 \leq b$ are constants independent of the ball. This property is often referred to as the `power law' for measures and ensures that the measure of any given ball does not depend too much on where the ball is centred. The condition also ensures that the measure of any given ball behaves in the way we would expect when it is scaled up or down. As a result of the power law alone, it is easily verified that if $m$ is $\delta$-Ahlfors regular
then
\begin{equation}
	\label{regcon1} \hdim \Omega = \delta.
\end{equation}

To define the second desirable property of the measures in which we are interested, we must first describe the collection of affine hyperplanes in projective space. The property from the classical setting that we will mimic is the idea that our measure must not be concentrated `too near' certain generic hyperplanes in $\R^{N+1}$. If this were not the case, then the measure might `live' too near to the rational points which we are trying to avoid.

Obviously, the projective space $\Proj^N(\R)$ is not a vector space, but we may define (affine) hyperplanes inside the space in the following way. We say that $\mathcal{P}$ is a \textit{projective hyperplane} if it can be expressed in the form $\mathcal{P} = \phi(\mathcal{L}\setminus\left\{\mathbf{0}\right\})$ for some $N$-dimensional (vector) hyperplane $\mathcal{L}$ of $\R^{N+1}$. Indeed, a projective hyperplane is then simply a subspace of projective space of codimension $1$, as one would expect. For example, when $N=1$ a projective hyperplane takes the form of a single point in $\Proj^1(\R)$, and when $N=2$ one can visualise a projective hyperplane as a great circle on the surface of the unit $2$-sphere with antipodal points associated; that is, a copy of $\Proj^1(\R)$ embedded inside $\Proj^2(\R)$.

Now, let $\mathcal{P}$ denote a generic projective hyperplane in $\Proj^N(\R)$, and for any $\epsilon > 0$ let $\mathcal{P}^{(\epsilon )}$ denote the $\epsilon$-neighbourhood of $\mathcal{P}$. To be
precise, 
\[
\mathcal{P}^{(\epsilon )}:=\left\{\xp \in \Proj^N(\R):
\inf_{\yp \, \in \, \mathcal{P}}\delta_\infty(\xp, \yp) < \epsilon \right\}.
\]
We say a measure $m$ is \emph{projectively $(c, \eta)$-absolutely decaying} if there exist strictly  positive constants $c, \eta$ and $r_0$ such that for any such projective hyperplane $\mathcal{P}$, any $\epsilon >0$, $\xp \in \Omega$ and $r< r_0$ we have
\begin{equation*}
	m(\hat{B}(\xp, r) \cap \mathcal{P}^{(\epsilon )}) \leq c \left(
	\frac{\epsilon}{r} \right)^{\eta}
	m(\hat{B}(\xp, r)).
\end{equation*}
We say $m$ is \textit{projectively absolutely decaying} if it
is projectively $(c, \eta)$-absolutely decaying for some $c, \eta >0$.

In the one-dimensional situation (i.e., $\Omega \subset \Proj^1(\R)$),
one can check that the Ahlfors regular property implies the projectively absolutely decaying property. However, this is not true in general.

\subsection{Example - the measure $\hat\mu_\infty$}
\label{sec:powerproj}

One can verify that $\hat{\mu}_\infty$ satisfies the power law with exponent $N$ and $r_0\leq\sin1$ and is projectively absolutely decaying with $\eta=N$. We will explicitly need the first of these properties in our proofs so for completion we include the details here.

  For any $\yp \in \PNR$ define $\pi(\yp) \subset \RNI$ to be the $N$-dimensional linear subspace of $\RNI$ orthogonal to the line $\phi^{-1}(\yp)\cup\left\{\zero\right\}$. Now, the subspace $\pi(\yp)$ naturally divides $\RNI$ into two open `half-spaces', each with boundary $\pi(\yp)$. Let us arbitrarily denote by $\Pi_1(\yp)$ the half-space containing the point $(1, 0, \ldots, 0)$,  or, if $(1, 0, \ldots, 0)\in \pi(\yp)$, then denote by $\Pi_1(\yp)$ the half-space containing $(0, 1, 0, \ldots, 0)$, or, if $\pi(\yp)$ contains both $(1, 0, \ldots, 0)$ and $(0, 1, 0, \ldots, 0)$, then denote by $\Pi_1(\yp)$ the half-space containing $(0, 0, 1, 0, \ldots, 0)$, and so on. The other half-space shall be denoted $\Pi_2(\yp)$.

For any ball $\hat{B}(\xp, r)$ of radius $r \leq \sin1$ in $\PNR$ the $(N+1)$-dimensional Lebesgue measure $\lambda$ of $\phi^{-1}(\hat{B}(\xp, r)) \cap B(\zero, 1)$ is then given by twice the volume of the hyperspherical cone $C_1$ given by $$C_1 \, = \, (\phi^{-1}(\hat{B}(\xp, r)) \cap \Pi_1(\xp))\cup\left\{\zero\right\}.$$ In particular (see for example \cite{Li}), we have
$$
\lambda(\phi^{-1}(\hat{B}(\xp, r)) \cap B(\zero, 1)) \, = \, 2 \cdot \lambda(B(\zero, 1)) \cdot\frac{2 \, \Gamma((N+1)/2)}{\sqrt{\pi}\, \Gamma(N/2)} \int_0^\theta \sin^{N-1}z \, dz,
$$
where $\theta \in (0, 1)$ is such that $r=\sin\theta$ and $\Gamma$ is the gamma function. Since for $z \in (0, 1)$ the inequalities
$$
z \, \geq \, \sin z \, \geq \, z - \frac{z^3}{6} \, \geq \frac{5z}{6}
$$
hold it follows that for $r < r_0 = \sin 1$ we have
$$
\frac{1}{N}\left(\frac{6}{5}\right)^{N} r^N \, \geq \, \frac{1}{N}\theta^N \, \geq \, \int_0^\theta \sin^{N-1}z \, dz \, \geq \,  \frac{1}{N}\left(\frac{5}{6}\right)^{N-1} \theta^N  \, \geq \, \frac{1}{N}\left(\frac{5}{6}\right)^{N-1} r^N.
$$
Thus, by formula (\ref{eqn:projmeasure}) we have
$$
a r^N \, \leq \, \hat{\mu}(\hat{B}(\xp, r)) \, \leq  b r^N,
$$
where
$$
a \, = \, \frac{4 \cdot 5^{N-1} \, \Gamma((N+1)/2)}{N\sqrt{\pi} \, 6^{N-1}\, \Gamma(N/2)} \quad \text{ and } \quad b \, = \, \frac{4 \cdot 6^{N} \, \Gamma((N+1)/2)}{N\sqrt{\pi} \, 5^{N}\, \Gamma(N/2)}.
$$




\section{Preliminaries for the proof of  Theorem \ref{thm:omegadim}}
\label{sec:framework}

The proof of Theorem \ref{thm:omegadim} makes use of the general
framework developed in \cite{MR2231044} for establishing dimension
statements for a large class of badly approximable sets. In this
section we provide a simplification of the framework that is
geared towards the particular application we have in mind. In
turn, this will avoid excessive referencing to the conditions
imposed in \cite{MR2231044} and thereby improve the clarity of our
exposition.

Let $(X, d)$ be a metric
space and $(Y, d)$ be a compact subspace of $X$ which
supports a non-atomic finite measure $m$. Let $\mathcal{R} :=
\left\{ R_{a} \in X : a \in J \right\} $ be a family of
subsets $R_{a}$ of $X$ indexed by an infinite countable set
$J$.  The sets $R_{a}$ will be referred to as the
\emph{resonant sets}. Next, let $\beta :J\rightarrow
\mathbb{R}_{>0}:a \mapsto \beta_{a}$ be a positive
function on $J$ such that the number of $a \in J$ with $
\beta_{a} $ bounded above is finite. Thus, $ \beta_{a} $
tends to infinity as $a$ runs through $J$. Finally, let $\rho: \R^+ \rightarrow\R^+:t \rightarrow \rho(t)$ be a decreasing function such that $\rho(t) \rightarrow 0$ as $t \rightarrow \infty$. Assume that for $t > 1$ sufficiently large and any integer $n \geq 1$ we have
\begin{equation}\label{eqn:rho}
\ell_1(t) \, \leq \, \frac{\rho(t^{n})}{\rho(t^{n+1})} \, \leq \,  \ell_2(t),
\end{equation}
where $\ell_1$ and $\ell_2$ are lower and upper bounds depending only on $t$ such that $\ell_1 \rightarrow \infty$ as $t \rightarrow \infty$.

 We are now in the
position to define the badly approximable set. Let
\begin{equation*}
\Bad(\mathcal{R}, \beta, \rho) := \left\{x \in Y : \exists \text{
} c(x)>0 \text{ s.t. } d(x, R_{a} ) \geq
c(x)\rho(\beta_{a}) \text {  } \forall \text{ } a \in
J \right\}  \ ,
\end{equation*}
where   $d(x, R_{a} ):= \inf_{r \in R_{a}}d(x, r)$.
Loosely speaking, $\Bad(\mathcal{R}, \beta, \rho)$ consists of points in
$Y $ that `stay clear' of the family $ \mathcal{R} $ of
resonant sets by a factor governed by $\rho(\beta)$.

The framework described in \cite{MR2231044} aims to determine the conditions under which $\hdim
\Bad(\mathcal{R}, \beta, \rho)  =  \hdim Y $.
 To specify these conditions we first require to introduce some notation.
For any fixed integer $t>1$ and any integer  $n \geq 1$, let $B_n
:= \left\{ x \in Y : d(c, x)
\leq \rho(t^n) \right\}$ denote a generic closed ball in $Y$ of radius $\rho(t^n)$
with centre $c$ in $Y$.  For any $\theta \in
\mathbb{R}_{>0}$, let $\theta B_n := \left\{ x \in Y : d(c,
x)
\leq \theta\rho(t^n) \right\}$ denote the ball $B_n$ scaled by $\theta$.
Finally,  let $J(n) := \left\{ a \in J : t^{n-1} \leq
\beta_{a} < t^n \right\}$. The following statement is a
simple consequence of  combining Theorem 1 and Lemma 7 of
\cite{MR2231044} and realises the above mentioned goal.

\begin{thmktv}
	\label{thm:1fromthat}
	Let $(X, d)$ be
	a metric space and $(Y, d)$ be a compact subspace of $X$
	which supports of a $\delta$-Ahlfors
	regular measure $m$ and let $\beta$ and $\rho$ satisfy the conditions stipulated above. Then, for $t$ sufficiently large, any
	$\theta \in \mathbb{R}_{>0}$, any $ n \geq 1 $ and any ball $B_n $
	there exists a collection $\mathcal{C}(\theta B_n)$ of
	disjoint balls $2 \theta B_{n+1}$ contained within $\theta B_n$ such
	that $$ \# \mathcal{C}(\theta B_n) \geq \kappa_1 \,  \left(\frac{\rho(t^n)}{\rho(t^{n+1})}\right)^{\delta}, $$
	where $\kappa_1>0$ is an absolute constant independent of $t$ and $n$.
	In addition, suppose for  some  $\theta \in \mathbb{R}_{>0}$ we also have that
	\begin{equation}
	\label{eq:cond2}
	\# \left\{ 2 \theta B_{n+1} \subset \mathcal{C}(\theta B_n):
	\min_{a \in J(n+1)} d(c, R_{a}) \leq 2
	\theta \rho(t^{(n+1)}) \right\} \leq \kappa_2 \left(\frac{\rho(t^n)}{\rho(t^{n+1})}\right)^{\delta}  \, ,
	\end{equation}
	where $\kappa_2$  is an absolute constant independent of $t$ and $n$ satisfying $0< \kappa_2 < \kappa_1$.  Furthermore, suppose
	\begin{equation}
	\label{eq:cond3}
	\hdim \, \left( \cup_{a \in J} R_{a} \right) < \delta  \, .
	\end{equation}
	Then
	\begin{equation*}
	\hdim \, \Bad(\mathcal{R}, \beta, \rho) = \delta  \, .
	\end{equation*}
\end{thmktv}

\section{Proof of Theorem \ref{thm:omegadim}}

\subsection{The strategy}

With reference to \S\ref{sec:framework}, we may set
\begin{eqnarray*}
	& X:=\PNR \, ,\quad Y :=\Omega \, ,\quad  d:= \delta_\infty\ , \quad J:=\left\{\mathbf{q} \in
	\mathbb{Z}^{N+1}\setminus\left\{\mathbf{0}\right\}  \right\}, \quad \rho(t):= t^{-(N+1)/N}, \\
	& a:=\mathbf{q} \in J, \quad
	R_a:= \qp=\phi(\mathbf{q}) \in \PNQ, \quad m:=\hat{\mu} \quad
	\rm{and} \quad    \beta_a:=H(\qp).  \
\end{eqnarray*}
It then follows that
$$\Bad(\mathcal{R}, \beta, \rho) = \Bad_\Proj^\infty(N) \cap \Omega, $$
and so the proof of Theorem \ref{thm:omegadim} is reduced to  showing that the conditions of
Theorem~KTV are satisfied.

For $t > 1$ and $n \geq 1$, let  $\hat{B}_n $ be a generic closed ball
of radius $\rho(t^n)=t^{-n(N+1)/N}$ and centre in $\Omega$, and so $\rho$ satisfies condition (\ref{eqn:rho}).  Assuming the measure $m$ is $\delta$-Ahlfors regular, for $t $ sufficiently large
and any $\theta \in \mathbb{R}_{>0}$ Theorem~\ref{thm:1fromthat} immediately provides a collection $\mathcal{C}(\theta \hat{B}_n)$ of disjoint
balls $2 \theta \hat{B}_{n+1}$ contained within $\theta \hat{B}_n$ such that
\begin{equation*}
\label{eq:cond1}
\# \mathcal{C}(\theta \hat{B}_n) \geq \kappa_1 \, \left(\frac{\rho(t^n)}{\rho(t^{n+1})}\right)^{\delta}\: = \: \kappa_1 t^{\delta(N+1)/N} \,  .
\end{equation*}
 We now endeavor to show that the additional condition
(\ref{eq:cond2}) on the collection $\mathcal{C}(\theta \hat{B}_n)$ is
satisfied.

We are required to find an upper bound for the number of possible disjoint balls $2\theta \hat{B}_{n+1}$, all lying inside a given fixed ball $\theta\hat{B_n}$, that can contain a projective rational $\qp \in \Omega\cap \theta\hat{B_n}$ with $t^n < H(\qp) \leq t^{n+1}$. To calculate this upper bound we will need to prove a lemma.

\subsection{A Symplex Lemma}
\label{sec:symplex}

In order to apply the framework outlined in the previous section we must first provide an analogue of the `symplex lemma' from classical Diophantine approximation.

We are given a ball $\hat{B}(\xp, r)$ of suitably small radius in $\PNR$ and $N+1$ distinct projective rationals $\qIp, \ldots, \qNIp$ in $\hat{B}(\xp, r)$ with $t^n < H(\qip) \leq t^{n+1}$ for $1 \leq i \leq N+1$. We wish to construct the analogue of an $N$-simplex within projective space with `vertices' at the points $\qIp, \ldots, \qNIp$. To do this we need to describe a method designed to mimic the process of taking the Euclidean convex hull.

With reference to \S\ref{sec:powerproj}, for each given projective rational $\qip \in \hat{B}(\xp, r)$ let $\qie \in \RNI$ be the unique point contained in the intersection $\phi^{-1}(\qip)\cap \Pi_1(\xp) \cap S^N$, where $S^N$ is the standard unit $N$-sphere in $\R^{N+1}$. We refer to $\qie$ as \textit{the canonical representative of $\qip$ in $\phi^{-1}(\hat{B}(\xp, r))$}. To be precise, assume the projective rational $\qip$ has homogeneous coordinates given by $\qip= [p_1^{(i)}, \ldots, p_{N+1}^{(i)}]$. Then, denoting by $\qie^{+}$ and $\qie^{-}$ the two points of intersection of $\phi^{-1}(\qip)$ with $S^N$; i.e.,
$$
\qie^{+}: \, = \, \frac{(p_1^{(i)}, \ldots, p_{N+1}^{(i)})}{\norm{(p_1^{(i)}, \ldots, p_{N+1}^{(i)})}}, \quad \text{ and } \quad \qie^{-}: \, = \, \frac{-(p_1^{(i)}, \ldots, p_{N+1}^{(i)})}{\norm{(p_1^{(i)}, \ldots, p_{N+1}^{(i)})}};
$$
then we define
$$
\qie \, = \, \begin{cases}
\qie^{+}, & \text{ if } \qie^{+} \in \Pi_1(\xp). \\
\qie^{-}, & \text{ if } \qie^{-} \in \Pi_1(\xp). \
\end{cases}
$$

Next, consider the Euclidean $(N+1)$-symplex $C_0:=C_{0}(\qIp, \ldots, \qNIp)$ in $\RNI$ formed by taking the convex hull of the $N+1$ canonical representatives $\qIe, \ldots, \qNIe$ and the origin.  For notational convenience let $C=C_{0} \setminus \left\{\zero\right\}$. We then define the \textit{projective $N$-symplex} $S:=S(\qIp, \ldots, \qNIp)$ \textit{with vertices} $\qIp, \ldots, \qNIp$ to be the image under the quotient map of the set $C$, i.e., $S=\phi(C)$. Note that $C$ is $\lambda$-measurable, since it is a convex subset of $\RNI$, 
and so $S$ is in turn $\hat{\mu}_\infty$-measurable. The reason for using the set $C$ as a base for the construction is the following. By design and the convexity of $\RNI$ we have
$$S \subset \hat{B}(\xp, r)  \quad \text{ and } \quad C\subset \phi^{-1}(S) \, \cap \, B(\zero, 1) \, \cap\, \Pi_1(\xp).$$
In particular, this implies
\begin{equation}
\label{eqn:simplex1}
\hat{\mu}_\infty(\hat{B}(\xp, r)) \, \geq \, \hat{\mu}(S)
\end{equation}
and
\begin{equation}
\label{eqn:simplex2}
\hat{\mu}_\infty(S) \, \geq \, \frac{2\lambda(C)}{\lambda(B(\zero, 1))}.
\end{equation}
Furthermore, $\lambda(C)=0$ if and only if the representatives $\qIe, \ldots, \qNIe$ lie in some $N$-dimensional linear subspace of $\RNI$, which occurs if and only if the projective rationals $\qIp, \ldots, \qNIp$ lie on an ($(N-1)$-dimensional) projective hyperplane.

Assume now that $\qIp, \ldots, \qNIp$ do not lie on some projective hyperplane, and so the volume of $C_0$ is given by
$$
\text{Vol}(C_0) \, = \, \infabs{\frac{\det\left[\qIe, \ldots, \qNIe \right]}{(N+1)!}} \, > \, 0.
$$
It follows that
\begin{eqnarray*}
	\lambda(C) & \geq & \frac{1}{(N+1)!} \left( \prod_{i=1}^{N+1} \, \norm{(p_1^{(i)}, \ldots, p_{N+1}^{(i)})} \right)^{-1} \\	
	& = & \frac{1}{(N+1)!} \left( \prod_{i=1}^{N+1} \, \sum_{j=1}^{N+1} \left(p_j^{(i)}\right)^2 \right)^{-1/2} \\
	& \geq & \frac{1}{(N+1)! (N+1)^{(N+1)/2}} \prod_{i=1}^{N+1} H(\qip)^{-1}. \	
\end{eqnarray*}

\subsection{Completion of the proof}\label{sec:done}

Returning to our proof, assume we have $N+1$ distinct projective rationals $\qIp, \ldots, \qNIp$ in $\Omega\cap \theta\hat{B_n}$ for which
\begin{equation}\label{eqn:kadicbound}
t^n < H(\qip) \leq t^{n+1} \quad \rm{ for } \, 1 \leq i \leq N+1,
\end{equation}
and that $\qIp, \ldots, \qNIp$ do not all lie on some projective hyperplane in $\PNR$.  Then, for sufficiently large $n$ the projective symplex lemma tells us that the projective symplex in $\PNR$ subtended by the points $\qIp, \ldots, \qNIp$ has $\hat{\mu}_\infty$-measure at most $c_N\prod_{i=1}^{N+1} H(\qip)^{-1}$, where $c_N$ is some strictly positive absolute constant depending only on $N$. In turn, the bounds given by (\ref{eqn:kadicbound}) imply that this quantity is bounded below by $c_N \,t^{-(n+1)(N+1)}$.

Furthermore, as $\hat{\mu}_\infty$ is $N$-Ahlfors regular (see \S\ref{sec:powerproj}) then $\hat{\mu}_\infty(\theta\hat{B_n})$ is bounded above by $b_N \, \theta\rho( t^n)^N= b_N \,\theta t^{-n(N+1)}$, for some strictly positive absolute constant $b_N$. Therefore, if we take
$$
\theta:= \frac{c_N}{2\, b_N \, t^{N+1}},
$$
we have a contradiction since the projective symplex subtended by $\qIp, \ldots, \qNIp$ must lie in the ball $\theta\hat{B_n}$ by local convexity. As a result, the projective rationals $\qIp, \ldots, \qNIp$ must all lie on some projective hyperplane $\mathcal{P}=\mathcal{P}(n)$ passing through $\theta \hat{B}_n$.

Therefore, a necessary condition for a particular ball $2\theta \hat{B}_{n+1}$ to contain a projective rational in $\Omega$ satisfying (\ref{eqn:kadicbound}) is that the ball intersects $\mathcal{P}$. Now, fix
$$
\epsilon \: = \: \epsilon(n): \: = \: \frac{4\theta}{t^{(n+1)(N+1)/N}}.
$$
Then the number of possible disjoint balls $2\theta \hat{B}_{n+1}$ inside $\theta\hat{B_n}$ that can contain a projective rational $\qp \in \PNR \cap\Omega$ satisfying (\ref{eqn:kadicbound}) is bounded above by the cardinality of the set
$$
\left\{2\theta \hat{B}_{n+1} \subset \mathcal{C}(\theta \hat{B}_n): \: 2\theta \hat{B}_{n+1} \subset \mathcal{P}^{(\epsilon)} \right\}.
$$
Since the balls in question are disjoint and the measure $\hat{\mu}$ on $\Omega$ is assumed to be $\delta$-Ahlfors regular (for some constants $a$ and $b$ with $0 < a \leq 1 \leq b$) and projectively $(c, \eta)$-absolutely decaying, we have that for sufficiently large $n$ that this cardinality is bounded above by
\begin{eqnarray*}
	\frac{\hat{\mu}(2\theta \hat{B}_{n+1} \cap \mathcal{P}^{(\epsilon)})}{\hat{\mu}(2\theta \hat{B}_{n+1} )} & \leq &
	c\left(\frac{\epsilon}{\theta t^{-n(N+1)/N}}\right)^\eta \hat{\mu}(\theta \hat{B}_n)/ \hat{\mu}(2\theta \hat{B}_{n+1} )\\
	& = & 2^{-\delta}a^{-1}bc \left(4t^{-(N+1)/N}\right)^\eta t^{\delta(N+1)/N} \\
	& = & 2^{2\eta-\delta}a^{-1}bc \: t^{(\delta-\eta)(N+1)/N} \\
	& < & \kappa_2\, t^{\delta(N+1)/N}, \
\end{eqnarray*}
when $t$ is sufficiently large, for any strictly positive $\kappa_2<\kappa_1$. Thus, our collection $\mathcal{C}(\theta \hat{B}_n)$ satisfies the condition (\ref{eq:cond2}).

Finally, in the context of the KTV framework observe that our resonant sets $\mathcal{R}$ are simply points in $\PNQ$. Hence, we trivially have that $\hdim \, \left( \cup_{a \in J} R_{a} \right)=0$ and so (\ref{eq:cond3}) also holds, meaning all of the conditions of Theorem KTV are satisfied and our proof is complete.

\section{Proof of Theorem $\ref{thm:Barak}$}

Let $\psi: \N \rightarrow \R_{>0}$ be a non-increasing arithmetic function for which
	\begin{equation} \label{eqn:barakconv}
	\sum_{m=1}^\infty \, m^{\frac{N+1}N\eta-1} \psi(m)^{\eta} \: < \: \infty.
	\end{equation}
For notational convenience let
$$
\mathcal{W}_\Omega^N(\psi): \, = \, \left\{\xp \in \Omega: \: \delta_\infty(\xp, \qp)  \leq \psi(H(\qp)) \, \text{ for infinitely many } \qp \in \PNQ   \right\}.
$$
Then, Theorem $\ref{thm:Barak}$ is precisely the statement that $\hat{\mu}(\mathcal{W}_\Omega^N(\psi))=0$. Since $\Omega$ is compact and the measure of a countable union of null sets is itself a null set we may restrict our attention to a generic ball $\hat{B}_{(0)}$ in $\PNR$ of radius strictly less than $ \sin 1$. It suffices to show that $\hat{\mu}(\mathcal{W}_\Omega^N(\psi) \cap \hat B_{(0)})=0$ for any such ball $\hat{B}_{(0)}$ for which this intersection is non-empty.  In turn, by the Borel-Cantelli lemma it suffices to show that
\begin{equation}
\label{eqn:BCconv}
\hat{\mu} \left( \bigcup_{\qp \, \in \, \PNQ \, \cap \, \hat B_{(0)} \, \cap \, \Omega} \: \hat B (\qp, \, \psi(H(\qp)))   \right) \: < \: \infty.
\end{equation}
In particular, we will show that for some $n_0 \geq 1$ and some $t>1$ sufficiently large we have
$$
\hat{\mu}\left(   W_n  \right)   \: \ll \:   \left(t^{n\frac{N+1}N} \psi(t^n)\right)^{\eta}
$$
for all $n \geq n_0$, where
$$
W_n: \: = \: \bigcup_{\substack{\qp \, \in \, \PNQ \, \cap \, \hat B_{(0)} \, \cap \, \Omega \\ t^n \, \leq \, H(\qp) \, < \, t^{n+1}}}  \: \hat B (\qp, \, \psi(H(\qp)));
$$
for, then by Cauchy condensation the left hand side of (\ref{eqn:BCconv}) is bounded above by
\begin{equation*}
\label{eqn:endgame}
\sum_{n=n_0}^\infty \, \hat{\mu}(W_n) \: \ll \: \sum_{n=n_0}^\infty \, \left(t^{n\frac{N+1}N} \psi(t^n)\right)^{\eta} \: \ll \: \sum_{m=1}^\infty \, m^{\frac{N+1}N\eta-1} \psi(m)^{\eta} \, < \, \infty.
\end{equation*}

Firstly, since $\hat{B}_{(0)}\cap \Omega$ is compact, for each $n \geq n_0$ it follows from Vitali's famous covering lemma that there exists a finite disjoint collection $\mathcal{B}_{(n)}$ of balls $\hat B_{(n)}$, with centres in $\hat{B}_{(0)}\cap \Omega$ and each of radius
$$
r_n \, = \, \frac{c_N}{4(1+2t)\, b_N} \, t^{-(n+1)\frac{(N+1)}{N}},
$$
for which
$$
\hat{B}_{(0)}\cap \Omega \quad \subset \quad \bigcup_{\hat B_{(n)} \, \in \, \mathcal{B}_{(n)}} \: (1+2t)\hat B_{(n)}.
$$
Here, $b_N$ and $c_N$ are constants as defined in \S\ref{sec:done}.
Now, assume that there is a projective  rational $\qp$ satisfying
\begin{equation}\label{eqn:diadic}
t^n \, \leq \, H(\qp) \, < \, t^{n+1}.
\end{equation}
and
$$
\hat B (\qp, \, \psi(H(\qp))) \, \cap \,  (1+2t)\hat B_{(n)} \, \neq \, \emptyset.
$$
Notice that we for any constant  $\gamma>0$ we may choose $n_0$ sufficiently large so that for $n \geq n_0$ we have
$$
\psi(t^n) < \gamma t^{-n\frac{N+1}{N}}.
$$
If this were not the case we would in view of (\ref{eqn:BCconv}) have a contradiction. Specifically, we may choose $n_0$ large enough so that $B (\qp, \, \psi(H(\qp)))$ is contained in $ 2(1+2t)\hat B_{(n)}$, and certainly so large that $\qp$ is contained in  $ 2(1+2t)\hat B_{(n)}$. However, as in the proof of Theorem \ref{thm:omegadim},  the projective symplex lemma then implies that all projective rationals $\qp \in B_0$ which satisfy (\ref{eqn:diadic}) and are contained in $2(1+2t)\hat B_{(n)}$ must lie on some projective hyperplane $\mathcal{P}_n$ passing through $2(1+2t)\hat B_{(n)}$.

Finally, set $\epsilon:= \psi(t^n)$. Then since the measure $\hat{\mu}$ is projectively absolutely $\eta$-decaying measure $m$ we have for $n \geq n_0$ that
\begin{eqnarray*}
	\hat{\mu}(W_n)  & = & \hat{\mu} \left( W_n \cap  \bigcup_{\hat B_{(n)}  \in  \mathcal{B}_{(n)}}  2(1+2t)\hat B_{(n)}     \right) \\ & \leq & \sum_{\hat B_{(n)}  \in  \mathcal{B}_{(n)}}
	\hat{\mu} \left(  \mathcal{P}_n^{(\epsilon)}  \cap   2(1+2t)\hat B_{(n)}   \right) \\
	&  \leq &
	c \left(\frac{4(1+2t)\, b_N}{c_N} \, t^{(n+1)\frac{(N+1)}{N}} \,\psi(t^n) \right)^\eta     \sum_{\hat B_{(n)}  \in  \mathcal{B}_{(n)}} \,  2(1+2t)\hat B_{(n)}. \
\end{eqnarray*}
Since $\hat{\mu}$ satisfies the power law it follows that
$$
\sum_{\hat B_{(n)}  \in  \mathcal{B}_{(n)}} \,  2(1+2t)\hat B_{(n)} \, \ll \, m(\hat{B}_{(0)}\cap \Omega) \, \ll \, 1,
$$ for all $n \geq n_0$. This in turn implies that there is an absolute constant $c_0$ depending on upon $N$ and $t$ for which
$$
\hat{\mu}(W_n) \, \leq \, c_0 \, \left(t^{n\frac{N+1}N} \psi(t^n)\right)^{\eta},
$$
and our proof is complete.

\noindent\emph{\bf Acknowledgments.} We would like to thank Prof. Victor Beresnevich and Simon Kristensen for many useful discussions. In addition, the first author believes a paper is never complete without reserved thanks for Prof. Sanju Velani.

\def\cprime{$'$} \def\cprime{$'$} \def\cprime{$'$}

\end{document}